\documentclass[12pt,a4paper]{article}
\usepackage{}
\usepackage{indentfirst}
\usepackage{amsfonts}
\usepackage{amssymb}
\usepackage{mathrsfs}
\usepackage{amsmath}
\usepackage{amsthm}
\usepackage{enumerate}
\usepackage{cite}
\usepackage[all]{xy}
\allowdisplaybreaks
\newtheorem{defn}{Definition}[section]
\newtheorem{thm}[defn]{Theorem}
\newtheorem{lem}[defn]{Lemma}
\newtheorem{prop}[defn]{Proposition}
\newtheorem{cor}[defn]{Corollary}
\newtheorem{eg}[defn]{Example}
\newtheorem{re}[defn]{Remark}
\newcommand{\bdefn}{\begin{defn}}
\newcommand{\edefn}{\end{defn}}
\newcommand{\bthm}{\begin{thm}}
\newcommand{\ethm}{\end{thm}}
\newcommand{\blem}{\begin{lem}}
\newcommand{\elem}{\end{lem}}
\newcommand{\bprop}{\begin{prop}}
\newcommand{\eprop}{\end{prop}}
\newcommand{\bcor}{\begin{cor}}
\newcommand{\ecor}{\end{cor}}
\newcommand{\beg}{\begin{eg}}
\newcommand{\eeg}{\end{eg}}
\newcommand{\bre}{\begin{re}}
\newcommand{\ere}{\end{re}}
\newcommand{\bpf}{\begin{proof}}
\newcommand{\epf}{\end{proof}}
\newcommand{\F}{{\rm{\bf F}}}
\newcommand{\id}{{\rm id}}

\newcommand{\supercite}[1]{\textsuperscript{\cite{#1}}}
\newcommand{\benu}{\begin{enumerate}}
\newcommand{\eenu}{\end{enumerate}}
\newcommand{\bc}{\begin{center}}
\newcommand{\ec}{\end{center}}
\newcommand{\bea}{\begin{eqnarray}}
\newcommand{\eea}{\end{eqnarray}}
\newcommand{\Bea}{\begin{eqnarray*}}
\newcommand{\Eea}{\end{eqnarray*}}
\newcommand{\beq}{\begin{equation}}
\newcommand{\eeq}{\end{equation}}
\newcommand{\Beq}{\begin{equation*}}
\newcommand{\Eeq}{\end{equation*}}
\newcommand{\bspl}{\begin{split}}
\newcommand{\espl}{\end{split}}
\newcommand\relphantom[1]{\mathrel{\phantom{#1}}}

\numberwithin{equation}{section}

\begin{document}
\title{\bf The construction and deformation of Hom-Novikov superalgebras}
\author{\normalsize \bf Bing Sun, Liangyun Chen,  Yan Liu}
\date{\small{ School of Mathematics and Statistics, Northeast Normal University,\\ Changchun,  130024,  CHINA}}
\maketitle
\begin{abstract}
We study a twisted generalization of Novikov superalgebras, called Hom-Novikov superalgebras. It is shown
 that two classes of Hom-Novikov superalgebras can be constructed from Hom-supercommutative algebras
  together with derivations and Hom-Novikov superalgebras with Rota-Baxter operators, respectively. We show that quadratic Hom-Novikov superalgebras are Hom-associative superalgebras and the sub-adjacent Hom-Lie superalgebras of
  Hom-Novikov superalgebras are 2-step nilpotent. Moreover,
  we develop the 1-parameter formal deformation theory of Hom-Novikov superalgebras.
\bigskip

\end{abstract}
 \footnote[0]{Corresponding author(L. Chen):
chenly640@nenu.edu.cn.} \footnote[0]{Supported by  NNSF of China
(Nos. 11171055 and 11471090),  NSF of  Jilin province (No.
201115006). }
\section{Introduction}
 Novikov algebras were introduced in connection with Hamiltonian operators in the formal variational calculus and the
 Poisson brackets of hydrodynamic-type.
They were used to construct the Virasoro-type Lie algebras. When
Gel'fand and Dorfman\supercite{Gel'fand1,Gel'fand2,Gel'fand3}
studied the following operator:
     $$H_{ij}=\sum_{k}c_{ijk}u_{k}^{(1)}+d_{ijk}u_{k}^{(0)}\frac{\rm d}{{\rm d}x}, \quad c_{ijk}\in{ \bf {C}},
   \quad d_{ijk}=c_{ijk}+c_{jik},$$
  they  gave the definition of Novikov algebras.
In 1987,  Zel'manov\supercite{Zel'manov} began to study Novikov
algebras and
 proved that the dimension of finite-dimensional simple  Novikov algebras over a
field of  characteristic zero is one. In algebras, what are paid
attention to by mathematician are  classifications and structures,
but so far we haven't got the systematic theory for general Novikov
algebras. In 1992,   Osborn\supercite{Osborn1,Osborn2} had  finished
the classification of  infinite simple  Novikov algebras with
nilpotent elements over  a field of  characteristic zero and finite
simple Novikov algebras with nilpotent elements over a field of
 characteristic $p>0$. In 1995, Xu\supercite{Xu3,Xu4,Xu5,Xu6} developed his theory and
 got  the classification of simple Novikov algebras
 over an algebraically closed field of  characteristic zero.
  Bai and   Meng \supercite{Bai,Bai1,Bai2} did a series of researches
 on low dimensional Novikov algebras, such as the  structure  and
 classification. Chen construct  two kinds of  Novikov algebras\supercite{Chen1,Chen2}. Recently, people
 obtained some properties in Novikov superalgebras\supercite{Zhu,Ni,Kang}.

To be more precise, recall that a left-symmetric superalgebra is a
$\mathbb{Z}_2$-graded vector space A together with a binary
operation $\mu:A\times A\rightarrow A$ satisfying \beq
(xy)z-x(yz)=(-1)^{|x||y|}((yx)z-y(xz)),\label{eq1} \eeq for $x\in
A_{|x|}$, $y\in A_{|y|}$, $z\in A_{|z|}$, $|x|,|y|,|z|\in
\mathbb{Z}_2$.  Here and in what follows we often write
$\mu(x,y)=xy$. In the other words, if $(x,y,z)=(xy)z-x(yz)$ denotes
the associator, then (\ref{eq1}) shows that $(x,y,z)$ is
supersymmetric with respect to two left variables $x$ and $y$, hence
named left-symmetric superalgebra. A Novikov superalgebra is a
left-symmetric superalgebra $A$ that satisfies the additional
property \beq (xy)z=(-1)^{|y||z|}(xz)y, \,\,for \,x\in A, y\in
A_{|y|}, z\in A_{|z|}.\label{eq2} \eeq

A quadratic Novikov superalgebra, introduced in \cite{Ni}, is a Novikov superalgebra with a symmetric nondegenerate even invariant bilinear form.
The motivation for studying quadratic Novikov superalgebras comes from the fact that Lie
or associative algebras with symmetric nondegenerate even invariant bilinear forms have important applications in several areas of mathematics
and physics, such as the structure theory of finite-dimensional semi-simple Lie algebras,
the theory of complete integrable Hamiltonian systems and the classification of statistical
models over two-dimensional graphs.

 Yau in \cite{Yau} introduced Hom-Novikov algebras, in which the two defining identities are twisted by a linear map.
 It turned out that Hom-Novikov algebras can be constructed from Novikov algebras, commutative Hom-associative algebras
  and Hom-Lie algebras along with some suitable linear maps. And Yuan and You in \cite{Yuan} introduced quadratic  Hom-Novikov algebras.
  Later, Zhang, Hou and Bai in \cite{Zhang}
defined a Hom-Novikov superalgebra as a twisted generalization of Novikov superalgebras.

The deformation is a tool to study a mathematical object by deforming it into a family of the same kind of objects depending on a certain parameter. The deformation theory was introduced by Gerstenhaber for rings and algebras\supercite{Gerstenhaber1, Gerstenhaber2, Gerstenhaber3, Gerstenhaber4}, by Kubo and Taniguchi for Lie triple systems\supercite{Kubo&Taniguchi}, by Ma, Chen and Lin for Hom-Lie Yamaguti algebras\supercite{Ma1} and Hom-Lie triple systems\supercite{Ma2}, by Bai and Meng for Novikov algebras\supercite{Bai}. They studied 1-parameter formal deformations and established the connection between the cohomology groups and infinitesimal deformations.

The purpose of this paper is to consider the realization and the 1-parameter formal deformation theory of Hom-Novikov superalgebras based on some work in \cite{Bai, Ma2, Yau, Yuan}. The paper is organized as follows. Section 2 concerns the realization of Hom-Novikov superalgebras. It is shown
 that two classes of Hom-Novikov superalgebras can be constructed from Hom supercommutative algebras
  together with derivations and Hom-Novikov superalgebras with Rota-Baxter operators, respectively. In Section 3, we introduce
  the notion of quadratic Hom-Novikov superalgebras and the relationship between Hom-Novikov superalgebras and sub-adjacent Hom-Lie superalgebras.
   We also show quadratic Hom-Novikov superalgebras are Hom-associative superalgebras and the sub-adjacent Hom-Lie superalgebras of Hom-Novikov superalgebras are
   2-step nilpotent. Section 4 is dedicated to the 1-parameter formal deformation theory of Hom-Novikov superalgebras. Moreover, we define low orders coboundary
   operator and give low order cohomology groups of Hom-Novikov superalgebras.
 We show that the cohomology group is suitable for this 1-parameter formal deformation theory.

Throughout this paper $\F$ denotes an arbitrary field.

\section{Hom-Novikov superalgebra}
Let $(A,\cdot)$ be an algebra over field $\F$. $A$ is said to be a
superalgebra if the underlying vector space of $A$ is
$\mathbb{Z}_2$-graded (i.e., $A=A_{\bar{0}}\oplus A_{\bar{1}}$,
where $A_{\bar{0}}$ and $A_{\bar{1}}$ are vector subspaces of $A$)
and $A_{\alpha}\cdot A_{\beta}\subset A _{\alpha+\beta}$, $\forall
\alpha,\beta\in \mathbb{Z}_2$. An element $x\in A$ is called
homogeneous if $x\in A_{\bar{0}}\cup  A_{\bar{1}}$. In this work,
all elements are supposed to be homogeneous unless otherwise stated.
For a homogeneous
 element $x$ we shall use the standard notation $|x|\in \mathbb{Z}_2=\{\bar{0},\bar{1}\}$ to indicate its degree.

\bdefn{\rm\supercite{Ammar}}
A multiplicative Hom-Lie superalgebra is a triple $(A, [\cdot, \cdot]_A, \alpha)$
 consisting of a $\mathbb{Z}_2$-graded vector space $L$,  an even bilinear map $[\cdot, \cdot]:A\times A\rightarrow A$ and an even
 algebraic morphism $\alpha:A\rightarrow A$ satisfying
\beq
 \alpha([x, y])=[\alpha(x), \alpha(y)],
 \eeq
\beq
[x, y]=-(-1)^{|x||y|}[y, x],
\eeq
\beq
(-1)^{|x||z|}[\alpha(x), [y, z]]+(-1)^{|y||x|}[\alpha(y), [z, x]]+(-1)^{|z||y|}[\alpha(z), [x, y]]=0,
\eeq
where $x, y$ and $z$ are homogeneous elements in $A$.
 Lie superalgebras are examples of Hom-Lie superalgebras in which $\alpha$ is the identity map.
\edefn

\bdefn{\rm\supercite{Ammar}}
A Hom-associative superalgebra
 consists of a $\mathbb{Z}_2$-graded vector space $A$, a linear self-map $\alpha$ and
an even bilinear map $\mu :A\times A\rightarrow A$, satisfying
$$\alpha(xy)=\alpha(x)\alpha(y)\quad (multiplicativity)$$
and
$$(xy)\alpha(z)=\alpha(x)(yz)\quad (Hom\!-\!associativity),$$
for $x,y,z\in A$.
\edefn

\bdefn{\rm\supercite{Zhang}}\label{D1} A Hom-Novikov superalgebra is
a triple  $(A,\mu,\alpha)$ consisting of a $\mathbb{Z}_2$-graded
vector space $A$, an even bilinear map $\mu:A\times A\rightarrow A$
and an even homomorphism $\alpha:A\rightarrow A$ satisfying \beq
\alpha(xy)=\alpha(x)\alpha(y),\label{eq3} \eeq \beq (xy)
\alpha(z)-\alpha(x)(yz)=(-1)^{|x||y|}((yx)\alpha(z)-\alpha(y)(xz)),\label{eq4}
\eeq \beq (xy) \alpha(z)=(-1)^{|y||z|}(xz)\alpha(y).\label{eq5} \eeq
\edefn

We see that Novikov
superalgebras are examples of Hom-Novikov superalgebras in which
$\alpha$ is the identity map. For a Hom-Novikov superalgebra
$(A,\mu,\alpha)$, we call $\mu$ the Hom-Novikov product of $A$.
Comparing with Definition \ref{D1} with (\ref{eq1}) and (\ref{eq2}), if and only if (\ref{eq3}) and (\ref{eq4}) are satisfied, then we call
$(A,\mu,\alpha)$ a Hom-left-symmetric superalgebra. In particular, a
Hom-Novikov superalgebra is a Hom-left-symmetric superalgebra that
also satisfies (\ref{eq5}).
\bdefn
Let $(A,\mu,\alpha)$ be a  Hom-Novikov superalgebra, which is called
 \begin{enumerate}[\rm(i)]
\item regular if $\alpha$ is an algebra automorphism;
\item involutive if $\alpha$ is an involution, i.e., $\alpha^{2}={\rm id}$.
\end{enumerate}
\edefn
\blem{\rm\supercite{Zhang}}
Let $(A,\mu,\alpha)$ be a  Hom-Novikov superalgebra and $[\cdot,\cdot]:A\times A\rightarrow A$ be a binary operation on $A$ defined by
$$[x,y]=xy-(-1)^{|x||y|}yx\,\, \,\,\,\forall x,y\in A.$$
Then $HLie(A)=(A,[\cdot,\cdot],\alpha)$ is a Hom-Lie superalgebra with the same twisting map $\alpha$, which is called the sub-adjacent Hom-Lie superalgebra of $A$.
\elem
Zhang,Hou and Bai in \rm{\cite{Zhang}} gave a way to construct Hom-Novikov superalgebras, starting from a Novikov superalgebra and an algebra endomorphism. In the following, we provide a construction of  Novikov superalgebras from Hom-Novikov superalgebras along with an algebra automorphism.

\bthm\label{thm2}
If $(A,\mu,\alpha)$ is an involutive Hom-Novikov superalgebra, then $(A, \alpha \circ \mu,\alpha)$ is a Novikov superalgebra.
\ethm
\bpf
For convenience, we write $x\ast y=\alpha(xy)$, for all $x,y \in A$. Hence, it needs to show
\beq
(x\ast y) \ast z-x\ast (y\ast z)=(-1)^{|x||y|}((y\ast x)\ast z-y\ast(x\ast z)),\label{eq8}
\eeq
\Beq
(x\ast y) \ast z=(-1)^{|y||z|}(x\ast z)\ast y,
\Eeq
for all $x,y,z \in A$. Since $(A,\mu,\alpha)$ is an involutive Hom-Novikov superalgebra, we have
$$(x\ast y) \ast z=\alpha(\alpha(xy)z)=\alpha^{2}(xy)\alpha(z)=(xy)\alpha(z)=(-1)^{|y||z|}(xz)\alpha(y)=(-1)^{|y||z|}(x\ast z) \ast y.$$
Furthermore,
\begin{align*}
(x\ast y) \ast z-x\ast (y\ast z)=&\alpha(\alpha(xy)z)-\alpha(x\alpha(yz))\\
=&(xy)\alpha(z)-\alpha(x)(yz)\\
=&(-1)^{|x||y|}((yx)\alpha(z)-\alpha(y)(xz))\\
=&(-1)^{|x||y|}((y\ast x)\ast z-y\ast(x\ast z)),
\end{align*}
which proves Equation (\ref{eq8}) and the proposition.
\epf
\bthm
If $(A,\mu,\alpha)$ is a regular Hom-Novikov superalgebra, then $(A, [\cdot,\cdot]_{\alpha^{-1}}=\alpha^{-1}\circ [\cdot,\cdot])$ is a Lie superalgebra, where
$[x,y]=xy-(-1)^{|x||y|}yx$, for $x,y\in A$. In particular, if $\alpha$ is an involution, then $(A, [\cdot,\cdot]_{\alpha}=\alpha \circ [\cdot,\cdot])$ is a Lie superalgebra.
\ethm
\bpf
For any $x,y,z\in A$, we have
\begin{align*}
&(-1)^{|x||z|}[x,[y,z]_{\alpha^{-1}}]_{\alpha^{-1}}\\
=&(-1)^{|x||z|}\alpha^{-1}(x[y,z]_{\alpha^{-1}}-(-1)^{|x||y|+|x||z|}[y,z]_{\alpha^{-1}}x)\\
=&(-1)^{|x||z|}\alpha^{-1}(x\alpha^{-1}(yz-(-1)^{|y||z|}zy)-(-1)^{|x||y|+|x||z|}\alpha^{-1}(yz-(-1)^{|y||z|}zy)x)\\
=&\alpha^{\!-\!2}((\!-\!1)^{|x||z|}\alpha(x)(yz)\!-\!(\!-\!1)^{|x||z|\!+\!|y||z|}\alpha(x)(zy)\!-\!(\!-\!1)^{|x||y|}(yz)\alpha(x)\!+\!(\!-\!1)^{|x||y|+|y||z|}(zy)\alpha(x)).
\end{align*}
Similarly, we have
\begin{align*}
&(-1)^{|x||y|}[y,[z,x]_{\alpha^{-1}}]_{\alpha^{-1}}\\
=&\alpha^{\!-\!2}((\!-\!1)^{|x||y|}\alpha(y)(zx)\!-\!(\!-\!1)^{|x||y|\!+\!|x||z|}\alpha(y)(xz)\!-\!(\!-\!1)^{|y||z|}(zx)\alpha(y)\!+\!(\!-\!1)^{|y||z|+|x||z|}(xz)\alpha(y)),\\
&(-1)^{|y||z|}[z,[x,y]_{\alpha^{-1}}]_{\alpha^{-1}}\\
=&\alpha^{\!-\!2}((\!-\!1)^{|y||z|}\alpha(z)(xy)\!-\!(\!-\!1)^{|x||y|\!+\!|y||z|}\alpha(z)(yx)\!-\!(\!-\!1)^{|x||z|}(xy)\alpha(z)\!+\!(\!-\!1)^{|x||y|+|x||z|}(yx)\alpha(z)).
\end{align*}
Then it follows from Equation (\ref{eq4}) that
$$(-1)^{|x||z|}[x,[y,z]_{\alpha^{-1}}]_{\alpha^{-1}}+(-1)^{|x||y|}[y,[z,x]_{\alpha^{-1}}]_{\alpha^{-1}}+(-1)^{|y||z|}[z,[x,y]_{\alpha^{-1}}]_{\alpha^{-1}}=0.$$
Clearly, $[x,y]=-(-1)^{|x||y|}[y,x]$. Since $\alpha^{-1}$ is an
automorphism, we have
$[x,y]_{\alpha^{-1}}=-(-1)^{|x||y|}[y,x]_{\alpha^{-1}}$, which
proves that $(A, [\cdot,\cdot]_{\alpha^{-1}})$ is a Lie
superalgebra. It follows immediately that $(A,
[\cdot,\cdot]_{\alpha})$ is also a Lie superalgebra since
$\alpha=\alpha^{-1}$ when $\alpha$ is an involution. \epf
\bthm\label{thm3} Let $(A,\mu,\alpha)$ be a  Hom-Novikov
superalgebra. Then $(A, \alpha \circ \mu,\alpha^{2})$ is a
Hom-Novikov superalgebra. \ethm \bpf For convenience, we write
$x\ast y=\alpha(xy)$, for all $x,y \in A$. Hence, it needs to show
\beq (x\ast y) \ast \alpha^{2} (z)-\alpha^{2} (x)\ast (y\ast
z)=(-1)^{|x||y|}((y\ast x)\ast \alpha^{2} (z)-\alpha^{2}
(y)\ast(x\ast z)),\label{eq12} \eeq \Beq (x\ast y) \ast \alpha^{2}
(z)=(-1)^{|y||z|}(x\ast z)\ast \alpha^{2} (y), \Eeq for all $x,y,z
\in A$. Since $(A,\mu,\alpha)$ is a Hom-Novikov superalgebra, we
have
$$(x\ast y) \ast \alpha^{2} (z)=\alpha^{2}((xy)\alpha(z))=(-1)^{|y||z|}\alpha^{2}((xz)\alpha(y))=(-1)^{|y||z|}(x\ast z) \ast \alpha^{2} (y).$$
Furthermore,
\begin{align*}
(x\ast y) \ast \alpha^{2} (z)-\alpha^{2} (x)\ast (y\ast z)=&\alpha^{2}((xy)\alpha(z)-\alpha(x)(yz))\\
=&(-1)^{|x||y|}\alpha^{2}((yx)\alpha(z)-\alpha(y)(xz))\\
=&(-1)^{|x||y|}((y\ast x) \ast \alpha^{2} (z)-\alpha^{2} (y)\ast (x\ast z)),
\end{align*}
which proves Equation (\ref{eq12}) and the proposition.
\epf

We now construct several classes of Hom-Novikov superalgebras. Yau in \cite{Yau} gave such construction of
Hom-Novikov algebras, starts with an  Hom-associative and commutative algebra $(A,\mu)$ and an even derivation $D:A \rightarrow A$. The new product
 \beq
 x\ast y=\mu(x,D(y))=xD(y),\label{eq6}
 \eeq
 for $x,y\in A$ makes $(A,*)$ into a Hom-Novikov algebra. To generalize this construction, we define a Hom supercommutative algebra to be a Hom-associative
 superalgebra whose multiplication is supercommutative. An even derivation on a Hom-associative superalgebra is defined in the usual way. Then we have the following result.

\bthm
Let $(A,\mu,\alpha)$ be a  Hom-supercommutative algebra and $D:A \rightarrow A$ be an even derivation such that $D\alpha=\alpha D$. Then $(A,*,\alpha)$ is
a Hom-Novikov superalgebra, where $\ast$ is defined as in (\ref{eq6}).
\ethm
\bpf
The multiplication of $\alpha$ with respect to $*$ in (\ref{eq6}) follows from the multiplication of $\alpha$ with respect to $\mu$ and the hypothesis
$D\alpha=\alpha D$. Next we check (\ref{eq4})
\begin{eqnarray*}
 &&(x*y)*\alpha(z)-\alpha(x)*(y*z)\\
 &=&(xD(y))D(\alpha(z))-\alpha(x)D(yD(z))\\
 &=&(xD(y))\alpha(D(z))-\alpha(x)(D(y)D(z))-\alpha(x)(yD^{2}(z))\\
 &=&-\alpha(x)(yD^{2}(z))\\
 &=&-(xy)\alpha(D^{2}(z)).
\end{eqnarray*}
The last  two equalities  follow from Hom-associativity. On the other hand,
\begin{eqnarray*}
 &&(-1)^{|x||y|}((y*x)*\alpha(z)-\alpha(y)*(x*z))\\
 &=&(-1)^{|x||y|}((yD(x))D(\alpha(z))-\alpha(y)D(xD(z)))\\
 &=&(-1)^{|x||y|}((yD(x))\alpha(D(z))-\alpha(y)(D(x)D(z))-\alpha(y)(xD^{2}(z)))\\
 &=&(-1)^{|x||y|}(\alpha(y)(D(x)D(z))-\alpha(y)(D(x)D(z))-\alpha(y)(xD^{2}(z)))\\
 &=&-(-1)^{|x||y|}\alpha(y)(xD^{2}(z))\\
 &=&-(-1)^{|x||y|}(yx)\alpha(D^{2}(z))\\
 &=&-(xy)\alpha(D^{2}(z)).
\end{eqnarray*}
Futhermore, we have
  \begin{eqnarray*}
  (x*y)*\alpha(z)&=&(xD(y))\alpha(D(z))=\alpha(x)(D(y)D(z))\\
  &=&(-1)^{|y||z|}\alpha(x)(D(z)D(y))\\
  &=&(-1)^{|y||z|}(xD(z))\alpha(D(y))\\
  &=&(-1)^{|y||z|}(x*z)*\alpha(y).
  \end{eqnarray*}
  Consequently, we prove the theorem.
\epf

\bthm
  Let $(A,\mu)$ be an associative
 superalgebra whose multiplication is supercommutative,  $\alpha:A \rightarrow A$ be an even algebra morphism, and $D:A \rightarrow A$ be
  an even derivation such that $D\alpha=\alpha D$. Then $(A,*,\alpha)$ is a Hom-Novikov superalgebra, where
\Beq
x\ast y=\alpha(xD(y)),\,\,\forall x,y\in A.
\Eeq
\ethm
\bpf
First, we have
 \begin{eqnarray*}
  \alpha(x*y)&=&\alpha\circ \alpha(xD(y))=\alpha(\alpha(x)D(\alpha(y)))=\alpha(x)*\alpha(y).
  \end{eqnarray*}
 On one hand, we have
  \begin{eqnarray*}
 &&(x*y)*\alpha(z)-\alpha(x)*(y*z)\\
 &=&\alpha(xD(y))*\alpha(z)-\alpha(x)*\alpha(yD(z))\\
 &=&\alpha((xD(y))*z-x*(yD(z)))\\
  &=&\alpha(\alpha((xD(y))D(z))-\alpha(xD(yD(z))))\\
 &=&\alpha^{2}((xD(y))D(z)-xD(yD(z)))\\
 &=&\alpha^{2}(x(D(y)D(z))-x(D(y)D(z))-x(yD^{2}(z)))\\
 &=&-\alpha^{2}(x(yD^{2}(z)))\\
 &=&-\alpha^{2}((xy)D^{2}(z)).
\end{eqnarray*}
Similarly, we have
  \begin{eqnarray*}
   &&(-1)^{|x||y|}((y*x)*\alpha(z)-\alpha(y)*(x*z))\\
 &=&(-1)^{|x||y|}\alpha((yD(x))*z-y*(xD(z)))\\
 &=&(-1)^{|x||y|}\alpha^{2}((yD(x))D(z)-yD(xD(z)))\\
 &=&(-1)^{|x||y|}\alpha^{2}(y(D(x)D(z))-y(D(x)D(z))-y(xD^{2}(z)))\\
 &=&-(-1)^{|x||y|}\alpha^{2}((yx)D^{2}(z))\\
 &=&-\alpha^{2}((xy)D^{2}(z)).
\end{eqnarray*}
Futhermore, we have
 \begin{eqnarray*}
  (x*y)*\alpha(z)&=&\alpha(xD(y))*\alpha(z)=\alpha^{2}((xD(y))D(z))\\
  &=&\alpha^{2}(x(D(y)D(z)))\\
  &=&(-1)^{|y||z|}\alpha^{2}(x(D(z)D(y)))\\
  &=&(-1)^{|y||z|}\alpha^{2}((xD(z)D(y))\\
  &=&(-1)^{|y||z|}(x*z)*\alpha(y).
  \end{eqnarray*}
 Consequently, we prove the theorem.
\epf

\bdefn
Let $(A,\ast,\alpha)$ be a Hom-superalgebra and let $\lambda\in {\F}$. If a linear map $P:A\rightarrow A$ satisfies
\begin{align*}
P(x)\ast P(y)=P(P(x)\ast y+x\ast P(y)+\lambda x\ast y), \,\,\,\forall x,y\in A,
\end{align*}
then $P$ is called a Rota-Baxter operator of weight $\lambda$ and $(A,\ast,\alpha, P)$ is called a Rota-Baxter Hom-superalgebra of weight $\lambda$.
\edefn
\bthm
Let $(A,*,\alpha,P)$ be a Rota-Baxter Hom-Novikov superalgebra of weight $\lambda$ and $P$ an even linear map. Assume that $\alpha$ and $P$ commute.
  Then $(A,\circ,\alpha,P)$ is a Hom-Novikov superalgebra, where the multiplication $\circ$ is defined as
\Beq
  x\circ y:=P(x)*y+x*P(y)+\lambda x*y,\,\,\,\forall x,y\in A.
\Eeq
\bpf
The multiplication of $\alpha$ with respect to $\circ$ follows from the multiplication of $\alpha$ with respect to $\ast$ and the hypothesis
$P\alpha=\alpha P$. For any $x,y,z\in A$, we have,
\begin{align*}
&(x\circ y)\circ \alpha(z)-\alpha(x)\circ (y\circ z)\\
=&P(P(x)\ast y+x\ast P(y)+\lambda x*y)\ast \alpha(z)+(P(x)\ast y+x\ast P(y)+\lambda x*y)\ast P(\alpha(z))\\
&+\lambda (P(x)\ast y+x\ast P(y)+\lambda x*y)\ast \alpha(z)-P(\alpha(x))\ast (P(y)\ast z+y\ast P(z)+\lambda y*z)\\
&-\alpha(x)\ast P(P(y)\ast z+y\ast P(z)+\lambda y*z)-\lambda \alpha(x)\ast (P(y)\ast z+y\ast P(z)+\lambda y*z)\\
=&(P(x)\ast P(y))\ast \alpha(z)+(P(x)\ast y)\ast \alpha(P(z))+(x\ast P(y))\ast \alpha(P(z))+\lambda (x*y)\ast \alpha(P(z))\\
&+\lambda (P(x)\ast y)\ast \alpha(z)+\lambda (x\ast P(y))\ast \alpha(z)+\lambda^{2}(x*y)\ast \alpha(z)-\alpha(P(x))\ast(P(y)\ast z)\\
&-\alpha(P(x))\ast (y\ast P(z))-\lambda \alpha(P(x))\ast (y\ast z)-\alpha(x)\ast (P(y)\ast P(z))-\lambda \alpha(x)\ast (P(y)\ast z)\\
&-\lambda \alpha(x)\ast (y\ast P(z))-\lambda^{2}\alpha (x)\ast (y\ast z)\\
=&(P(x)\ast P(y))\ast \alpha(z)\!-\!\alpha(P(x))\ast(P(y)\ast z)\!+\!(P(x)\ast y)\ast \alpha(P(z))\!-\!\alpha(P(x))\ast (y\ast P(z))\\
&+(x\ast P(y))\ast \alpha(P(z))\!-\!\alpha(x)\ast (P(y)\ast P(z))\!+\!\lambda (x*y)\ast \alpha(P(z))\!-\!\lambda \alpha(x)\ast (y\ast P(z))\\
&+\lambda (P(x)\ast y)\ast \alpha(z)-\lambda \alpha(P(x))\ast (y\ast z)+\lambda (x\ast P(y))\ast \alpha(z)-\lambda \alpha(x)\ast (P(y)\ast z)\\
&+\lambda^{2}(x*y)\ast \alpha(z)-\lambda^{2}\alpha (x)\ast (y\ast z)\\
=&(-1)^{|x||y|}\Big((P(y)\ast P(x))\ast \alpha(z)-\alpha(P(y))\ast(P(x)\ast z)+(y\ast P(x))\ast \alpha(P(z))\\
&\!-\!\alpha(y)\ast (P(x)\ast P(z))\!+\!(P(y)\ast x)\ast \alpha(P(z))\!-\!\alpha(P(y))\ast (x\ast P(z))\!+\!\lambda (y*x)\ast \alpha(P(z))\\
&-\lambda \alpha(y)\ast (x\ast P(z))+\lambda (y\ast P(x))\ast \alpha(z)-\lambda \alpha(y)\ast (P(x)\ast z)+\lambda (P(y)\ast x)\ast \alpha(z)\\
&-\lambda \alpha(P(y))\ast (x\ast z)+\lambda^{2}(y*x)\ast \alpha(z)-\lambda^{2}\alpha (y)\ast (x\ast z)\Big).
\end{align*}
Similarly, we have
 \begin{align*}
&(-1)^{|x||y|}\Big((y\circ x)\circ \alpha(z)-\alpha(y)\circ (x\circ z)\Big)\\
=&(-1)^{|x||y|}\Big((P(y)\ast P(x))\ast \alpha(z)-\alpha(P(y))\ast(P(x)\ast z)+(P(y)\ast x)\ast \alpha(P(z))\\
&\!-\!\alpha(P(y))\ast (x\ast P(z))\!+\!(y\ast P(x))\ast \alpha(P(z))\!-\!\alpha(y)\ast (P(x)\ast P(z))\!+\!\lambda (y*x)\ast \alpha(P(z))\\
&-\lambda \alpha(y)\ast (x\ast P(z))+\lambda (P(y)\ast x)\ast \alpha(z)-\lambda \alpha(P(y))\ast (x\ast z)+\lambda (y\ast P(x))\ast \alpha(z)\\
&-\lambda \alpha(y)\ast (P(x)\ast z)+\lambda^{2}(y*x)\ast \alpha(z)-\lambda^{2}\alpha (y)\ast (x\ast z)\Big).
 \end{align*}
Furthermore, on the one hand, we have
\begin{align*}
&(x\circ y)\circ \alpha(z)\\
=&P(P(x)\ast y+x\ast P(y)+\lambda x*y)\ast \alpha(z)+(P(x)\ast y+x\ast P(y)+\lambda x*y)\ast P(\alpha(z))\\
&+\lambda (P(x)\ast y+x\ast P(y)+\lambda x*y)\ast \alpha(z)\\
=&(P(x)\ast P(y))\ast \alpha(z)\!+\!(P(x)\ast y)\ast \alpha(P(z))\!+\!(x\ast P(y))\ast \alpha(P(z))\!+\!\lambda (x*y)\ast \alpha(P(z))\\
&+\lambda (P(x)\ast y)\ast \alpha(z)+\lambda (x\ast P(y))\ast \alpha(z)+\lambda^{2}(x*y)\ast \alpha(z)\\
=&(-1)^{|y||z|}\Big((P(x)\ast z)\ast \alpha(P(y))+(P(x)\ast P(z))\ast \alpha(y)+(x\ast P(z))\ast \alpha(P(y))\\
&+\lambda (x\ast P(z))\ast \alpha(y)+\lambda (P(x)\ast z)\ast \alpha(y)+\lambda (x\ast z)\ast \alpha(P(y))+\lambda^{2}(x*z)\ast \alpha(y)\Big).
\end{align*}
On the other hand, we have
\begin{align*}
&(-1)^{|y||z|}\Big((x\circ z)\circ \alpha(y)\Big)\\
=&(-1)^{|y||z|}\Big((P(x)\ast P(z))\ast \alpha(y)+(P(x)\ast z)\ast \alpha(P(y))+(x\ast P(z))\ast \alpha(P(y))\\
&+\lambda (x*z)\ast \alpha(P(y))+\lambda (P(x)\ast z)\ast \alpha(y)+\lambda (x\ast P(z))\ast \alpha(y)+\lambda^{2}(x*z)\ast \alpha(y)\Big).
\end{align*}
Hence, the conclusion holds.
\epf
\ethm

\section{Quadratic Hom-Novikov superalgebras}
\bdefn{\rm\supercite{Liu}}
 Let $A$ be a Hom-Lie superalgebra. A bilinear form $B$ on $A$ is
 \begin{enumerate}[\rm(i)]
\item supersymmetric if $B(x,y)=(-1)^{|x||y|}B(y,x)$, $\forall x, y\in A$;
\item nondegenerate if $x,y\in A$ satisfies $B(x,y)=0$, $\forall  y\in A$, then $x=0$;
\item invariant if $B([x,y],z)=B(x,[y,z])$, $\forall x,y,z\in A$;
\item even if $B(A_{\bar{0}},A_{\bar{1}})=B(A_{\bar{1}},A_{\bar{0}})=\{0\}$.
\end{enumerate}
 \edefn
\bdefn{\rm\supercite{Liu}}
A quadratic Hom-Lie superalgebra is a quadruple $(A,[\cdot,\cdot],\alpha, B)$ such that  $(A,[\cdot,\cdot],\alpha)$ is a Hom-Lie superalgebra  with an even, supersymmetric, non-degenerate and invariant bilinear form $B$ on $A$ satisfying
\beq
B(\alpha(x),y)=B(x,\alpha(y)),\,\,\forall x,y\in A.\label{eq10}
\eeq

We recover quadratic Lie superalgebras when $\alpha=\rm{id}.$
\edefn
\bdefn
A quadratic Hom-Novikov superalgebra is a quadruple $(A,\mu,\alpha, B)$ such that  $(A,\mu,\alpha)$ is a Hom-Novikov superalgebra  with an even, supersymmetric, non-degenerate bilinear form $B$ on $A$ satisfying
\beq\label{eq11}
B(\alpha(x),yz)=B(xy,\alpha(z)),\,\,\forall x,y,z\in A.
\eeq

We recover quadratic Novikov superalgebras when $\alpha=\rm{id}.$
\edefn
\bthm\label{thm1}
Let $(A,\mu,\alpha, B)$ be a quadratic Hom-Novikov superalgebra and $HLie(A)=(A,[\cdot,\cdot],\alpha)$ the sub-adjacent Hom-Lie superalgebra of $A$. If $\alpha$ is an automorphism satisfying
\beq
B(\alpha(x),y)=B(x,\alpha(y)),\,\,\forall x,y\in A,\label{eq9}
\eeq
then $(A,[\cdot,\cdot],\alpha, B_{\alpha})$ is a quadratic Hom-Lie superalgebra, where $B_{\alpha}(x,y)=B(\alpha(x),y).$
\ethm
\bpf
Since $B$ is a nondegenerate even bilinear form and $\alpha$ is an automorphism, $B_{\alpha}$ is a nondegenerate even
bilinear form on A. For all $x, y, z\in A$, using the properties of $B$, we have
\begin{align*}
B_{\alpha}([x,y],z)=&B(\alpha([x,y]),z)=B([x,y],\alpha(z))\\
=&B(xy,\alpha(z))-(-1)^{|x||y|}B(yx,\alpha(z))\\
=&B(\alpha(x),yz)-(-1)^{|x||y|}(-1)^{|x||y|+|y||z|}B(\alpha(x),zy)\\
=&B(\alpha(x),[y,z])=B_{\alpha}(x,[y,z]).
\end{align*}
Hence $B_{\alpha}$ is invariant. Using supersymmetry of $B$ and Equation (\ref{eq9}), we have
$$B_{\alpha}(x,y)=B(\alpha(x),y)=(-1)^{|x||y|}B(y,\alpha(x))=(-1)^{|x||y|}B(\alpha(y),x)=(-1)^{|x||y|}B_{\alpha}(y,x),$$
which proves $B_{\alpha}$ is supersymmetry. Using Equation (\ref{eq9}) again, we have
$$B_{\alpha}(\alpha(x),y)=B(\alpha(\alpha(x)),y)=B(\alpha(x),\alpha(y))=B_{\alpha}(x,\alpha(y)),$$
which completes the proof.
\epf
\bcor
Let $(A,\mu, B)$ be a quadratic Novikov superalgebra with an automorphism satisfying Equation (\ref{eq9})  and $(A,[\cdot,\cdot])$  be the sub-adjacent Lie superalgebra. Then $(A,[\cdot,\cdot]_{\alpha}=\alpha \circ [\cdot,\cdot] ,\alpha, B_{\alpha})$ forms
a quadratic Hom-Lie superalgebra, where $B_{\alpha}(x,y)=B(\alpha(x),y).$
\ecor
\bpf
It is obvious that $(A, [\cdot,\cdot]_{\alpha}, \alpha)$ is a Hom-Lie superalgebra. Using the similar
arguments as those in the proof of Theorem \ref{thm1}, we get $B_{\alpha}$ is an even symmetric nondegenerate bilinear
form with Equation (\ref{eq10}) satisfied. It remains to show that $B_{\alpha}$ is invariant. For all $x, y, z\in A$,
using invariance and supersymmetry of $B$, we have
\begin{align*}
B_{\alpha}([x,y]_{\alpha},z)=&B(\alpha([x,y]_{\alpha}),z)=B([x,y]_{\alpha},\alpha(z))\\
=&B(\alpha(x)\alpha(y),\alpha(z))-(-1)^{|x||y|}B(\alpha(y)\alpha(x),\alpha(z))\\
=&B(\alpha(x),\alpha(y)\alpha(z))-(-1)^{|x||y|}(-1)^{|x||y|+|y||z|}B(\alpha(x),\alpha(z)\alpha(y))\\
=&B(\alpha(x),[y,z]_{\alpha})=B_{\alpha}(x,[y,z]_{\alpha}),
\end{align*}
which proves the invariance of $B_{\alpha}$ and the result.
\epf

\bthm
Let $(A,\mu,\alpha,B)$ be a quadratic Hom-Novikov superalgebra, where $\alpha$ is an involution
satisfying Equation (\ref{eq9}). Then $(A,\alpha\circ \mu, B)$ is a quadratic Novikov superalgebra.
\ethm
\bpf
$(A,\alpha \circ\mu)$ is a Novikov superalgebra by Theorem \ref{thm2}. It suffices to show that $B$
is invariant under the operation $\alpha\circ \mu$. For all $x, y, z\in A$, we have
$$B(x, \alpha(y)\alpha(z))=B(\alpha(x), yz)=B(xy,\alpha(z))=B(\alpha(x)\alpha(y), z),$$
which completes the proof.
\epf

\bthm
Let $(A,\mu,\alpha, B)$ is a quadratic Hom-Novikov superalgebra, where $\alpha$ is an automorphism
satisfying Equation (\ref{eq9}). Then $(A,\ast=\alpha\circ \mu,\alpha^{2},B_{\alpha^{2}})$ is a quadratic Hom-Novikov superalgebra, where $B_{\alpha^{2}}(x,y)=B(\alpha^{2}(x),y).$
\ethm
\bpf
It follows from Theorem \ref{thm3} that $(A,\ast,\alpha^{2})$ forms a Hom-Novikov superalgebra. Since
$B$ is a nondegenerate even bilinear form on $A$ and $\alpha$ is an automorphism, $B_{\alpha^{2}}$ is a nondegenerate even bilinear
form. For all $x, y, z\in A$, using the hypothesis, we have
$$B_{\alpha^{2}}(x,y)=B(\alpha^{2}(x),y)=B(x,\alpha^{2}(y))=(-1)^{|x||y|}B(\alpha^{2}(y),x)=(-1)^{|x||y|}B_{\alpha^{2}}(y,x),$$
thus, $B_{\alpha^{2}}$ is supersymmetric. Moreover, we have
\begin{align*}
B_{\alpha^{2}}(\alpha^{2}(x),y\ast z)=&B(\alpha^{4}(x),\alpha(y)\alpha(z))=B(\alpha^{3}(x),\alpha^{2}(y)\alpha^{2}(z))=B(\alpha(x)\alpha(y),\alpha^{4}(z))\\
=&(-1)^{|x||z|+|y||z|}B(\alpha^{4}(z),x\ast y)=B_{\alpha^{2}}(x\ast y,\alpha^{2}(z)),
\end{align*}
which proves the invariance of $B_{\alpha^{2}}$ and the theorem.
\epf
\bcor
Let $(A,\mu,\alpha, B)$ be a quadratic Hom-Novikov superalgebra, where $\alpha$ is an automorphism
satisfying Equation (\ref{eq9}). Then $(A,\alpha^{n}\circ \mu, \alpha^{n},B_{\alpha^{n}})$ is a quadratic Hom-Novikov superalgebra, where $B_{\alpha^{n}}(x,y)=B(\alpha^{n}(x),y),$ $\forall n> 0$.
\ecor

Let $(A,\mu,\alpha)$ be a Hom-Novikov superalgebra, whose center is denoted by $Z(A)$ and defined by
$$Z(A)=\{x\in A|xy=yx=0, \,\,\forall y\in A\}.$$
Let $(A,[\cdot,\cdot],\beta)$ be a Hom-Lie superalgebra. The lower central series of $A$ is defined as usual, i.e., $A^{0}=A,\,A^{i}=[A,A^{i-1}],\,\forall i\geq 1$.
We call $A$ is $i$-step nilpotent if $A^{i}=0$ and $A^{i-1}\neq 0$. The center of the Hom-Lie
superalgebra is denoted by $C(A)$ and defined by
$$C(A)=\{x\in A|[x, y]=0,\,\, \forall y\in A\}.$$
\bthm\label{thm4}
Let $(A,\mu,\alpha, B)$ be a quadratic Hom-Novikov superalgebra. If  $\alpha$ is an automorphism, then $(A,\mu,\alpha)$ is a Hom-associative superalgebra.
\ethm
\bpf
Define $(x,y,z)=\alpha(x)yz-(xy)\alpha(z)$. For any $x,y,z,d\in A$, we have
\begin{align*}
B((x,y,z),\alpha(d))=&B(\alpha(x)(yz),\alpha(d))-B((xy)\alpha(z),\alpha(d))\\
=&B(\alpha^{2}(x),(yz)d)-B(\alpha(x)\alpha(y),\alpha(z\alpha^{-1}(d)))\\
=&B(\alpha^{2}(x),(yz)d)-B(\alpha^{2}(x),\alpha(y)(z\alpha^{-1}(d)))\\
=&-B(\alpha^{2}(x),(y,z,\alpha^{-1}(d))).
\end{align*}
Thus, we have
\begin{align*}
B((x,y,z),\alpha(d))=&-B(\alpha^{2}(x),(y,z,\alpha^{-1}(d)))\\
=&(-1)^{1+|x|(|y|+|z|+|d|)+|y||z|}B((z,y,\alpha^{-1}(d))),\alpha^{2}(x))\\
=&(-1)^{|x|(|y|+|z|+|d|)+|y||z|}B(\alpha^{2}(z),(y,\alpha^{-1}(d),x))\\
=&(-1)^{|x|(|y|+|z|+|d|)+|y||z|+|z|(|y|+|d|+|x|)}B((y,\alpha^{-1}(d),x),\alpha^{2}(z))\\
=&(-1)^{|x|(|y|+|z|+|d|)+|y||z|+|z|(|y|+|d|+|x|)+|y||d|}B((\alpha^{-1}(d),y,x),\alpha^{2}(z))\\
=&(-1)^{1+|x|(|y|+|z|+|d|)+|y||z|+|z|(|y|+|d|+|x|)+|y||d|}B(\alpha(d),(y,x,z))\\
=&(-1)^{1\!+\!|x|(|y|\!+\!|z|\!+\!|d|)\!+\!|y||z|\!+\!|z|(|y|\!+\!|d|\!+\!|x|)\!+\!|y||d|\!+\!|d|(|x|\!+\!|y|\!+\!|z|)\!+\!|x||y|}B((x,y,z),\alpha(d))\\
=&-B((x,y,z),\alpha(d)).
\end{align*}
It follows that $(x,y,z)=0$ by the  non-degeneracy of $B$.
\epf
\bthm
Let $(A,\mu,\alpha, B)$ be a quadratic Hom-Novikov superalgebra and $HLie(A)$ be the sub-adjacent Hom-Lie superalgebra. If  $\alpha$ is an automorphism, then $[x,y]\subseteq Z(A),$ for any $x,y\in HLie(A)$. As a
consequence, $H\!Lie(A)$ is 2-step nilpotent.
\ethm
\bpf
For any  $x, y, z\in A$, by Theorem \ref{thm4} we have
\begin{align*}
\alpha(z)[x,y]=&\alpha(z)(xy)-(-1)^{|x||y|}\alpha(z)(yx)=(zx)\alpha(y)-(-1)^{|x||y|}(zy)\alpha(x)\\
=&(zx)\alpha(y)-(zx)\alpha(y)=0.
\end{align*}
Using Equation (\ref{eq11}), we have
\begin{align*}
B([x,y]z,\alpha(d))=&B(\alpha[x,y],zd)=(-1)^{(|z|+|d|)(|x|+|y|)}B(zd,\alpha([x,y]))\\
=&(-1)^{(|z|+|d|)(|x|+|y|)}B(\alpha(z),d[x,y])=0,
\end{align*}
which implies $[x,y]\in Z(A)$ since $\alpha$ is an automorphism and $B$ is nondegenerate.
Hence, we have
$[HLie(A),HLie(A)]\subseteq Z(A)$. Obviously, $Z(A)\subseteq C(HLie(A))$. Then it follows that $HLie(A)$ is
2-step nilpotent.
\epf
\bre
A 2-step nilpotent quadratic Hom-Lie superalgebra $A$ admits a quadratic Hom-Novikov superalgebra. It suffices
to define a bilinear product on $A$ by $xy=\frac{1}{2}[x,y]$.
\ere
\section{1-parameter formal deformations of Hom-Novikov superalgebras}
\bdefn
Let $(A, \ast, \alpha)$ is a Hom-Novikov superalgebra. If an $n$-linear map $f:\underbrace{A\times\cdots\times A}_{n~{\rm times}}\rightarrow A$ satisfies
\begin{gather}
\alpha(f(x_1,\cdots,x_{n}))=f(\alpha(x_1),\cdots,\alpha(x_{n})),\label{cochain1}
\end{gather}
then $f$ is called an $n$-Hom-cochain on $A$. Denote by $C_{\alpha}^n(A,A)$ the set of all $n$-Hom-cochains, $\forall n\geq1$.
\edefn

\bdefn
For $n=1,2$, the coboundary operator $\delta^n_{hom}: C_{\alpha}^n(A,A) \rightarrow C_{\alpha}^{n+1}(A,A)$ is defined as follows.
\begin{align*}
\delta_{hom}^1f(x_1,x_2)=&(-1)^{x_{1}f}x_1\ast f(x_2)+f(x_1)\ast x_2-f(x_1\ast x_2);\\
\delta_{hom}^2f(x_1,x_2,x_3)=&f(\alpha(x_1),x_2\ast x_3)-(-1)^{|x_1||x_2|}f(\alpha(x_2),x_1\ast x_3)\\
&+(-1)^{|x_1||x_2|}f(x_2\ast x_1,\alpha(x_3))+(-1)^{|x_{1}||f|}\alpha(x_1)\ast f(x_2,x_3)\\
&-(-1)^{|x_1||x_2|+|x_{2}||f|}\alpha(x_2)\ast f(x_1,x_3)+(-1)^{|x_1||x_2|}f(x_2,x_1)\ast \alpha(x_3)\\
&-(-1)^{|x_2||x_3|}f(x_1\ast x_3,\alpha(x_2))-(-1)^{|x_2||x_3|}f(x_1,x_3)\ast \alpha(x_2).
\end{align*}

\edefn

It is not difficult to verify that $\delta^1_{hom}f$, $\delta^2_{hom}f$ satisfies (\ref{cochain1}). Thus, the coboundary operator $\delta^n_{hom}$ is well-defined.
\bthm
The coboundary operator $\delta^1_{hom}$, $\delta^2_{hom}$ defined above satisfies $\delta^{2}_{hom}\delta^1_{hom}=0.$
\ethm
\bpf
Suppose that $f\in C_{\alpha}^{1}(A,A)$, we have
\begin{align}
&\delta^{2}_{hom}\delta^{1}_{hom}f(x_1,x_2,x_3)\notag\\
=&\delta^{1}_{hom}f(\alpha(x_1),x_2\ast x_3)-(-1)^{|x_1||x_2|}\delta^{1}_{hom}f(\alpha(x_2),x_1\ast x_3)\notag\\
&+(-1)^{|x_1||x_2|}\delta^{1}_{hom}f(x_2\ast x_1,\alpha(x_3))+(-1)^{|x_{1}||f|}\alpha(x_1)\ast \delta^{1}_{hom}f(x_2,x_3)\notag\\
&-(-1)^{|x_1||x_2|+|x_{2}||f|}\alpha(x_2)\ast \delta^{1}_{hom}f(x_1,x_3)+(-1)^{|x_1||x_2|}\delta^{1}_{hom}f(x_2,x_1)\ast \alpha(x_3)\notag\\
&-(-1)^{|x_2||x_3|}\delta^{1}_{hom}f(x_1\ast x_3,\alpha(x_2))-(-1)^{|x_2||x_3|}\delta^{1}_{hom}f(x_1,x_3)\ast \alpha(x_2)\notag\\
=&(-1)^{|x_{1}||f|}\alpha(x_1)\ast f(x_2\ast x_3)+f(\alpha(x_1))\ast (x_2\ast x_3)-f(\alpha(x_1)\ast (x_2\ast x_3))\notag\\
&-(-1)^{|x_1||x_2|+|x_{2}||f|}\alpha(x_2)\ast f(x_1\ast x_3)-(-1)^{|x_1||x_2|}f(\alpha(x_2))\ast (x_1\ast x_3)\notag\\
&+(-1)^{|x_1||x_2|}f(\alpha(x_2)\ast (x_1\ast x_3))+(-1)^{|x_1||x_2|+(|x_1|+|x_{2}|)|f|}(x_2\ast x_1)\ast f(\alpha(x_3)\notag\\
&+(-1)^{|x_1||x_2|}f(x_2\ast x_1)\ast \alpha(x_3)-(-1)^{|x_1||x_2|}f((x_2\ast x_1)\ast \alpha(x_3))\notag\\
&+(-1)^{(|x_1+x_{2}|)|f|}\alpha(x_1)\ast (x_2\ast f(x_3))+(-1)^{|x_{1}||f|}\alpha(x_1)\ast (f(x_2)\ast x_3)\notag\\
&-(-1)^{|x_{1}||f|}\alpha(x_1)\ast f(x_2\ast x_3)-(-1)^{|x_1||x_2|+(|x_1|+|x_{2}|)|f|}\alpha(x_2)\ast (x_1\ast f(x_3))\notag\\
&-(-1)^{|x_1||x_2|+|x_{2}||f|}\alpha(x_2)\ast (f(x_1)\ast x_3)+(-1)^{|x_1||x_2|+|x_{2}||f|}\alpha(x_2)\ast f(x_1\ast x_3)\notag\\
&+(-1)^{|x_1||x_2|+|x_{2}||f|}(x_2\ast f(x_1))\ast \alpha(x_3)+(-1)^{|x_1||x_2|}(f(x_2)\ast x_1)\ast \alpha(x_3)\notag\\
&-(-1)^{|x_1||x_2|}f(x_2\ast x_1)\ast \alpha(x_3)-(-1)^{|x_2||x_3|+(|x_1|+|x_{3}|)|f|}(x_1\ast x_3)\ast f(\alpha(x_2))\notag\\
&-(-1)^{|x_2||x_3|}f(x_1\ast x_3)\ast \alpha(x_2)+(-1)^{|x_2||x_3|}f((x_1\ast x_3)\ast \alpha(x_2))\notag\\
&-(-1)^{|x_2||x_3|+|x_1||f|}(x_1\ast f(x_3))\ast \alpha(x_2)-(-1)^{|x_2||x_3|}(f(x_1)\ast x_3)\ast \alpha(x_2)\notag\\
&+(-1)^{|x_2||x_3|}f(x_1\ast x_3)\ast \alpha(x_2)\notag\\
=&f(\alpha(x_1))\ast (x_2\ast x_3)-(-1)^{|x_1||x_2|+|x_{2}||f|}\alpha(x_2)\ast (f(x_1)\ast x_3)\tag{a1}\\
&+(-1)^{|x_1||x_2|+|x_{2}||f|}(x_2\ast f(x_1))\ast \alpha(x_3)\tag{b1}\\
&-f(\alpha(x_1)\ast (x_2\ast x_3))+(-1)^{|x_1||x_2|}f(\alpha(x_2)\ast (x_1\ast x_3))\tag{a2}\\
&-(-1)^{|x_1||x_2|}f((x_2\ast x_1)\ast \alpha(x_3))\tag{b2}\\
&-(-1)^{|x_1||x_2|}f(\alpha(x_2))\ast (x_1\ast x_3)+(-1)^{|x_{1}||f|}\alpha(x_1)\ast (f(x_2)\ast x_3)\tag{a3}\\
&+(-1)^{|x_1||x_2|}(f(x_2)\ast x_1)\ast \alpha(x_3)\tag{b3}\\
&+(-1)^{|x_1||x_2|+(|x_1|+|x_{2}|)|f|}(x_2\ast x_1)\ast f(\alpha(x_3)+(-1)^{(|x_1+x_{2}|)|f|}\alpha(x_1)\ast (x_2\ast f(x_3))\tag{a4}\\
&-(-1)^{|x_1||x_2|+(|x_1|+|x_{2}|)|f|}\alpha(x_2)\ast (x_1\ast f(x_3))\tag{b4}\\
&-(-1)^{|x_2||x_3|+(|x_1|+|x_{3}|)|f|}(x_1\ast x_3)\ast f(\alpha(x_2)\tag{a5}\\
&+(-1)^{|x_2||x_3|}f((x_1\ast x_3)\ast \alpha(x_2))\tag{a6}\\
&-(-1)^{|x_2||x_3|+|x_1||f|}(x_1\ast f(x_3))\ast \alpha(x_2)\tag{a7}\\
&-(-1)^{|x_2||x_3|}(f(x_1)\ast x_3)\ast \alpha(x_2)\tag{a8}.
\end{align}
By (\ref{eq4}) and (\ref{eq5}), we have
\begin{align*}
&C1=(a1)+(b1)-(f(x_1)\ast x_2)\ast \alpha(x_3)=0\\
&C2=(a2)+(b2)+f((x_1\ast x_2)\ast \alpha(x_3))=0\\
&C3=(a3)+(b3)-(-1)^{|x_{1}||f|}(x_1\ast f(x_2))\ast \alpha(x_3)=0\\
&C4=(a4)+(b4)-(-1)^{(|x_1|+|x_{2}|)|f|}(x_1\ast x_2)\ast \alpha(f(x_3))=0\\
&C5=(a5)-(-1)^{|x_{1}||f|}(x_1\ast f(x_2))\ast \alpha(x_3)=0\\
&C6=(a6)-f((x_1\ast x_2)\ast \alpha(x_3))=0\\
&C7=(a7)+(-1)^{(|x_1|+|x_{2}|)|f|}(x_1\ast x_2)\ast \alpha(f(x_3))=0\\
&C8=(a8)+(f(x_1)\ast x_2)\ast \alpha(x_3)=0.
\end{align*}
 Thus,
\begin{align*}
&\delta^{2}_{hom}\delta^{1}_{hom}f(x_1,x_2,x_3)\\
=&(a1)+(b1)+(a2)+(b2)+(a3)+(b3)+(a4)+(b4)+(a5)+(a6)+(a7)+(a8)\\
=&C1+C2+C3+C4+C5+C6+C7+C8=0.
\end{align*}
The proof is completed.
\epf

For $n=1,2$, the map $f\in C_{\alpha}^n(A,A)$ is called an $n$-Hom-cocycle if $\delta^{n}_{hom}f=0$. We denote by $Z_{\alpha}^n(A,A)$ the subspace spanned by $n$-Hom-cocycles and $B_{\alpha}^n(A,A)=\delta^{n-1}_{hom}C_{\alpha}^{n-1}(A,A)$. Since $\delta^{2}_{hom}\delta^1_{hom}=0$, $B_{\alpha}^2(A,A)$ is a subspace of $Z_{\alpha}^2(A,A)$. Hence we can define a cohomology space $H_{\alpha}^2(A,A)$ of $(A, \ast, \alpha)$ as the factor space $Z_{\alpha}^2(A,A)/B_{\alpha }^2(A,A)$.

Let $(A, \ast, \alpha)$ be a Hom-Novikov superalgebra and $\F[[t]]$ be the ring of formal power series over $\F$. Suppose that $A[[t]]$ is the set of formal power series over $A$. Then for an $\F$-bilinear map $f:A\times A\rightarrow A$, it is natural to extend it to be an $\F[[t]]$-bilinear map $f:A[[t]]\times A[[t]]\rightarrow A[[t]]$ by
$$f\left(\sum_{i\geq0}x_it^i,\sum_{j\geq0}y_jt^j\right)=\sum_{i,j\geq0}f(x_i,y_j)t^{i+j}.$$
\bdefn
Let $(A, \ast, \alpha)$ be a Hom-Novikov superalgebra over $\F$. A 1-parameter formal deformation of $(A, \ast, \alpha)$ is a formal power series $g_t:A[[t]]\times A[[t]]\rightarrow A[[t]]$ of the form
$$g_t(x,y)=\sum_{i\geq0}G_i(x,y)t^i=G_0(x,y)+G_1(x,y)t+G_2(x,y)t^2+\cdots,$$
where each $G_i$ is an even $\F$-bilinear map $G_i:A\times A\rightarrow A$ $($extended to be $\F[[t]]$-bilinear$)$ and $G_0(x,y)=x\ast y$, such that the following identities hold
\begin{gather}
g_t(\alpha(x),\alpha(y))=\alpha\circ g_t(x,y),\label{deq1}\\
g_t(\alpha(x),g_t(y,z))-g_t(g_t(x,y),\alpha(z))=(-1)^{|x||y|}\{g_t(\alpha(y),g_t(x,z))-g_t(g_t(y,x),\alpha(z))\},\label{deq2}\\
g_t(g_t(x,y),\alpha(z))=(-1)^{|y||z|}g_t(g_t(x,z),\alpha(y)).\label{deq3}
\end{gather}
Conditions (\ref{deq1})-(\ref{deq3}) are called the deformation
equations of a Hom-Novikov superalgebra. \edefn

Note that $A[[t]]$ is a module over $\F[[t]]$ and $g_t$ defines the bilinear multiplication on $A[[t]]$ such that $A_t=(A[[t]],g_t,\alpha)$ is a Hom-Novikov superalgebra. Now we investigate the deformation Equations (\ref{deq1})-(\ref{deq3}).

Conditions (\ref{deq1})-(\ref{deq3}) are equivalent to the following
equations
\begin{gather}
G_i(\alpha(x),\alpha(y))=\alpha\circ G_i(x,y),\label{deq4}
\end{gather}
respectively, for $i=0,1,2,\cdots.$ The conditions (\ref{deq2}) and (\ref{deq3}) can be expressed as
\begin{align}
&\sum_{i,j\geq0}\Big(G_i(\alpha(x),G_j(y,z))-G_i(G_j(x,y),\alpha(z))-(-1)^{|x||y|}G_i(\alpha(y),G_j(x,z))\notag\\
&+(-1)^{|x||y|}G_i(G_j(y,x),\alpha(z))\Big)=0,\label{deq5}\\
&\sum_{i,j\geq0}\Big(G_i(G_j(x,y),\alpha(z))-(-1)^{|y||z|}G_i(G_j(x,z),\alpha(y))\Big)=0.\label{deq6}
\end{align}
For $n=0$, this means $A=A_{0}$ is a Hom-Novikov superalgebra.
For $n=1$, we obtain three equations for $G_{1}$:
\begin{align}
&G_i(\alpha(x),\alpha(y))=\alpha\circ G_i(x,y),\label{deq11}\\
&G_1(\alpha(x),y\ast z)\!-\!G_1(x\ast y,\alpha(z))\!-\!(-1)^{|x||y|}G_1(\alpha(y),x\ast z)\!+\!(-1)^{|x||y|}G_1(y\ast x,\alpha(z))+\notag\\
&\alpha(x)\!\ast\! G_1(y,z)\!\!-\!\!G_1(x,y)\!\ast\! \alpha(z)\!\!-\!\!(-\!1)^{|x||y|}\alpha(y)\!\ast\! G_1(x,z)\!\!+\!\!(-\!1)^{|x||y|}G_1(y,x)\!\ast\! \alpha(z)=0,\label{deq9}\\
&(G_1(x\ast y,\alpha(z))\!\!-\!\!(-1)^{|y||z|}G_1(x\ast z,\alpha(y))\!\!+\!\!G_1(x,y)\ast \alpha(z)\!\!-\!\!(-1)^{|y||z|}G_1(x,z)\ast \alpha(y).\label{deq10}
\end{align}
We call $G_{1}$ an infinitesimal deformation.
\beg\label{eg1}
A family of Hom-Novikov superalgebras $(A,\ast_{\xi},\alpha)$ is an infinitesimal deformation the Hom-Novikov superalgebra defined by Equation (\ref{eq6}) with a fixed derivation $D$, where $\ast_{\xi}$ is defined as
\beq
 x\ast_{\xi} y=xD(y)+\xi xy,\label{eq7}
 \eeq
 for all $x,y\in A,\,\xi\in \F.$ For Equation (\ref{eq7}), we let $G_{1}(x,y)=xy$, then Equation (\ref{deq9}) holds since
\begin{align*}
&\alpha (x) (y D(z))-(xD(y))\alpha (z)-(-1)^{|x||y|}\alpha(y)(xD(z))+(-1)^{|x||y|}(y D(x)) \alpha (z)\\
&+\alpha (x)D(yz)-(x y) D(\alpha (z))-(-1)^{|x||y|}\alpha(y)(D(x)z)+(-1)^{|x||y|}(y x) D(\alpha (z))=0
\end{align*}
and Equation (\ref{deq10}) holds since
\begin{align*}
&(x D(y)) \alpha (z)-(-1)^{|y||z|}(x D(z))\alpha(y)+(x y) D(\alpha (z))+(-1)^{|y||z|}(x z) D(\alpha (y))=0.
\end{align*}
Obviously, it is easy to verify that Equation (\ref{deq11}) holds.
\eeg
Hence, we have the following result.
\bthm
Let $(A,\mu,\alpha)$ be a  Hom-supercommutative algebra and $D:A \rightarrow A$ be an even derivation such that $D\alpha=\alpha D$. Then $(A,\ast_{\xi},\alpha)$ is
a Hom-Novikov superalgebra, where $\ast_{\xi}$ is defined as
\Beq
 x\ast_{\xi} y=xD(y)+\xi xy,
 \Eeq
 for all $x,y\in A,\,\xi\in \F.$
\ethm
For two $\F$-bilinear maps $f,g:A\times A\rightarrow A$ (extended to be $\F[[t]]$-bilinear), define a map $f\circ_{\alpha}g:A[[t]]\times A[[t]]\times A[[t]]\rightarrow A[[t]]$ by
\begin{align*}
f\circ_{\alpha}g(x,y,z)=&f(\alpha(x),g(y,z))-(-1)^{|x||y|}f(\alpha(y),g(x,z))+(-1)^{|x||y|}f(g(y,x),\alpha(z))\\
&-(-1)^{|y||z|}f(g(x,z),\alpha(y)).
\end{align*}
Using Equations (\ref{deq5}) and (\ref{deq6}), we have
$$\sum_{i+j=n}G_i\circ_{\alpha}G_j=0.$$

For $n=1$,
\beq
G_0\circ_{\alpha}G_1+G_1\circ_{\alpha}G_0=0.\label{deq7}
\eeq

For $n\geq2$,
\beq
-(G_0\circ_{\alpha}G_n+G_n\circ_{\alpha}G_0)=G_1\circ_{\alpha}G_{n-1}+G_2\circ_{\alpha}G_{n-2}+\cdots+G_{n-1}\circ_{\alpha}G_1.\label{deq8}
\eeq

By (\ref{deq4}) it follows that $G_i\in C_{\alpha}^2(A,A)$. It can also be verified that $G_i\circ_{\alpha}G_j\in C_{\alpha}^3(A,A)$. In general, if $f, g\in C_{\alpha}^2(A,A)$, then $f\circ_{\alpha}g\in C_{\alpha}^3(A,A)$. Note that the definition of coboundary operator, which implies $\delta_{hom}^2G_n=G_0\circ_{\alpha}G_n+G_n\circ_{\alpha}G_0,$ for $n=0,1,2\cdots.$ Hence (\ref{deq7}) and (\ref{deq8}) can be rewritten as
\begin{gather*}
\delta_{hom}^2G_1=0,\\
-\delta_{hom}^2G_n=G_1\circ_{\alpha}G_{n-1}+G_2\circ_{\alpha}G_{n-2}+\cdots+G_{n-1}\circ_{\alpha}G_1.
\end{gather*}
Then $G_1$ is a 2-Hom-cocycle.
\bdefn
Let $(A, \ast, \alpha)$ be a Hom-Novikov superalgebra. Suppose that $g_t(x,y)=\sum_{i\geq0}G_i(x,y)t^i$ and $g_t'(x,y)=\sum_{i\geq0}G_i'(x,y)t^i$ are two 1-parameter formal deformations of $(A, \ast, \alpha)$. They are called equivalent, denoted by $g_t\sim g_t'$, if there is a formal isomorphism of $\F[[t]]$-modules
$$\phi_t(x)=\sum_{i\geq0}\phi_i(x)t^i:(A[[t]],g_t,\alpha)\longrightarrow (A[[t]],g_t',\alpha),$$
where each $\phi_i:T\rightarrow T$ is an $\F$-linear map (extended to be $\F[[t]]$-linear) and $\phi_0=\id_A$, satisfying
\begin{gather*}
\phi_t\circ\alpha=\alpha\circ\phi_t,\\
\phi_t\circ g_t(x,y)=g_t'(\phi_t(x),\phi_t(y)).
\end{gather*}
When $G_1=G_2=\cdots=0$, $g_t=G_0$ is said to be the null deformation. A 1-parameter formal deformation $g_t$ is called trivial if $g_t\sim G_0$. A Hom-Novikov superalgebra $(T, [\cdot,\cdot,\cdot], \alpha)$ is called analytically rigid, if every 1-parameter formal deformation $g_t$ is trivial.
\edefn
\bthm
Let $g_t(x,y,z)=\sum_{i\geq0}G_i(x,y,z)t^i$ and $g_t'(x,y,z)=\sum_{i\geq0}G_i'(x,y,z)t^i$ be equivalent 1-parameter formal deformations of $(A, \ast, \alpha)$. Then $G_1$ and $G_1'$ belong to the same cohomology class in $H_{\alpha}^2(A,A)$.
\ethm
\bpf
Suppose that $\phi_t(x)=\sum_{i\geq0}\phi_i(x)t^i$ is the formal $\F[[t]]$-module isomorphism such that $\phi_t\circ\alpha=\alpha\circ\phi_t$ and
$$\sum_{i\geq0}\phi_i\left(\sum_{j\geq0} G_j(x,y)t^j\right)t^i =\sum_{i\geq0}G_i'\left(\sum_{k\geq0}\phi_k(x)t^k,\sum_{l\geq0}\phi_l(y)t^l \right)t^i.$$
It follows that
$$\sum_{i+j=n}\phi_i(G_j(x,y))t^{i+j}=\sum_{i+k+l=n}G_i'(\phi_k(x),\phi_l(y))t^{i+k+l}.$$
In particular,
$$\sum_{i+j=1}\phi_i(G_j(x,y))=\sum_{i+k+l=1}G_i'(\phi_k(x),\phi_l(y)),$$
that is,
$$G_1(x,y)+\phi_1(x\ast y)=\phi_1(x)\ast y+x\ast \phi_1(y)+G_1'(x,y).$$
Then $G_1-G_1'=\delta_{hom}^1\phi_1\in B_{\alpha}^2(A,A)$.
\epf

\bthm
Suppose that $(A, \ast, \alpha)$ is a Hom-Novikov superalgebra such that $H_{\alpha}^2(A,A)$ $= 0$. Then $(A, \ast, \alpha)$ is analytically rigid.
\ethm
\bpf
Let $g_t$ be a 1-parameter formal deformation of $(A, \ast, \alpha)$. Suppose that $g_t=G_0+\sum_{i\geq n}G_it^i$. Then
$$\delta_{hom}^2G_n=G_1\circ_{\alpha}G_{n-1}+G_2\circ_{\alpha}G_{n-2}+\cdots+G_{n-1}\circ_{\alpha}G_1=0,$$
that is, $G_n\in Z_{\alpha}^2(A,A)=B_{\alpha}^2(A,A)$. It follows that there exists $f_n\in C_{\alpha}^1(A,A)$ such that $G_n=\delta_{hom}^1f_n$.

Let $\phi_t=\id_A-f_nt^n:(A[[t]],g_t',\alpha)\longrightarrow (A[[t]],g_t,\alpha)$. Note that
$$\phi_t\circ \sum_{i\geq0}f_n^it^{in}=\sum_{i\geq0}f_n^it^{in}\circ \phi_t=\id_{A[[t]]}.$$
Then $\phi_t$ is a linear isomorphism. Moreover, $\phi_t\circ \alpha=\alpha\circ \phi_t$.

Now consider $g_t'(x,y)=\phi_t^{-1}g_t(\phi_t(x),\phi_t(y))$. It is straightforward to prove that $g_t'$ is a 1-parameter formal deformation of $(A, \ast, \alpha)$ and $g_t\sim g_t'$. Suppose that $g_t'=\sum_{i\geq 0}G_i't^i$. Then
$$(\id_T-f_nt^n)\left(\sum_{i\geq 0}G_i'(x,y)t^i\right)=\left(G_0+\sum_{i\geq n}G_it^i\right)(x-f_n(x)t^n,y-f_n(y)t^n),$$
i.e.,
\begin{align*}
&\sum_{i\geq 0}G_i'(x,y)t^i-\sum_{i\geq 0}f_n\circ G_i'(x,y)t^{i+n}\\
=&x\ast y-(f_n(x)\ast y+x\ast f_n(y))t^n+f_n(x)\ast f_n(y)t^{2n}\\
 &+\sum_{i\geq n}G_i(x,y)t^i-\sum_{i\geq n}(G_i(f_n(x),y)+G_i(x,f_n(y)))t^{i+n}+\sum_{i\geq n}G_i(f_n(x),f_n(y))t^{i+2n}.
\end{align*}
Then we have $G_1'=\cdots=G_{n-1}'=0$ and
$$G_n'(x,y)-f_n(x\ast y)=-(f_n(x)\ast y+x\ast f_n(y))+G_n(x,y).$$
Hence $G_n'=G_n-\delta_{hom}^1f_n=0$ and $g_t'=G_0+\sum_{i\geq n+1}G_i't^i$. By induction, this procedure ends with $g_t\sim G_0$, i.e., $(A, \ast, \alpha)$ is analytically rigid.
\epf

\end{document}